# Exploiting sparseness in damage characterization: A close look at the regularization techniques

Esmaeil Memarzadeh[1*], Dionisio Bernal[2], Martin D. Ulriksen[3]

[1,2] Civil and Environmental Engineering, Center for Digital Signal Processing, Northeastern University, Boston MA, USA

[3] Department of Energy Technology, Aalborg University, 6700 Esbjerg, Denmark

`memarzadeh@coe.neu.edu`

**Abstract.** The idea of exploiting sparseness in under-determined damage characterization problems is not new, and regularizations techniques that tend to promote sparseness, such as L1-norm minimization, have been investigated in the last ten years or so. Although various claims of merit have been made, two interconnected issues put these claims into question, and this paper brings some attention to the matter. The first is that the relationship between the structural parameters and the modal features previously considered has been linear and to ensure that the premise was closely realized, only very small damage severities have been considered. The second issue, intimately related to the first, is the fact that the noise, which has been typically taken as small relative to the "change in the features", is then unrealistically small. In problems where the damage is sufficiently large, the nonlinear dependence of the features on the parameters cannot be generally discarded. It is found that the attainable performance is much less "impressive" that what has been often claimed. The paper also examines the potential merit of using Lp-norm (0<p<1) minimization, instead of L1-norm minimization which, to the knowledge of the writers, has not been previously examined in damage characterization research. In this case we also find that, contrary to claims made in other areas, this norm does not lead to any general improvement over the performance attained by minimizing the L1-norm.

## 1 Introduction

Structural damage, which is a change in structural properties, causes changes in vibration characteristics such as the modal data. Sensitivity-based model updating methods can be used to characterize the changes in the structural properties by adjusting a model to minimize the discrepancy between the experimental data and the corresponding model predictions [1]. Modal parameters such as frequencies and mode shapes are commonly used as the features for damage characterization. Using frequencies is more common than mode shapes since frequencies can be estimated with more precision [2]. One major limitation of using the frequencies is that the number of frequencies is usually smaller than the number of parameters to be estimated/updated, which renders the problem under-determined. Therefore, the problem has an infinite set of solutions and regularization techniques are needed to be used. One idea that has been used recently is to promote sparsity. Damage usually occurs in a few members or sections of the structure. Therefore the small number of damage locations is the available prior knowledge that can be used as a criterion to look for solutions that have a smaller number of non-zeros. This is a well-known problem which is called a sparse recovery problem. A sparse vector has a relatively small number of non-zero elements. The number of non-zero elements in a vector is referred to as the L0-norm of that vector. L0-norm minimization looks for the solution with the least number of non-zeros among all the solutions. However, L0-norm minimization is non-convex, and L1-norm minimization, which is a convex problem, is used as a surrogate. L1-norm minimization has been used for damage characterization in the last ten years. Zhou et al. [2] first implemented the idea of using sparsity constraints in damage characterization by using L1-norm minimization. Hernandez [3] used sensitivity-based model updating method and applied L1-norm minimization using frequencies only. Hou [4] expanded this method to using the eigenvalues and the modes shapes.

In this study, we examine the use of regularization in damage characterization of under-determined problems. Following this introduction, section two briefly reviews the sensitivity method and highlights the fact that the relationship between the estimated modal features and the structural parameters is nonlinear. Section three looks at the regularization techniques. Section four contains a numerical example, and the paper closes with a summary in section five.

## 2 Sensitivity-based model updating

In this section, we briefly review the sensitivity method [1]. The modal features of a structure $f$ are a function of the structural parameters $\theta$. Writing the Taylor series expansion at an arbitrary set of parameters ($\theta \in R^{n \times 1}$), expanded about the nominal set of parameters ($\theta_0 \in R^{n \times 1}$), we get:

$$f(\theta) - f(\theta_0) = \frac{\delta f}{\delta \theta}(\theta - \theta_0) + \frac{\delta^2 f}{\delta \theta^2}\frac{(\theta - \theta_0)^2}{2} + \cdots \qquad (1)$$

where $f(\theta) \in R^{m \times 1}$ is the vector of modal features at a given set of structural parameters, and $\frac{\delta f}{\delta \theta} \in C^{m \times n}$ is the Jacobian at the nominal parameters. Taking $\theta - \theta_0$ as $x$, $f(\theta) - f(\theta_0)$ as $b$, $\frac{\delta f}{\delta \theta}$ as $A$, and the higher order terms as $hot$, we can re-write Eq. (1) as:

$$b = Ax + hot \qquad (2)$$

It is essential to mention that the relation between the features and parameter changes is nonlinear. In most previous studies on L1-norm minimization for damage characterization, the nonlinear terms were discarded by claiming that the relationship between small changes in features and small changes in parameters is linear [3]. However, this is not a valid assumption once the extent of damage is larger, and an iterative procedure is needed for a possible convergence to a solution.

## 3 Regularization

### 3.1 Sparsity

As it can be seen in Eq. (2), if the number of features is less than the number of parameters, that is, $m < n$, the problem will be under-determined, and it will have infinite solutions. We tend to promote sparsity and look for a solution that has the least number of non-zeros because damage usually happens in a few members or sections of the structure.

### 3.2 Minimum cardinality problem

The L0-norm, or the cardinality of a vector, is the number of its non-zero entries. Therefore, the following optimization problem will give the solution with the least number of non-zeros:

$$\min \|x\|_0 \ \ s.t \ Ax = b \qquad (3)$$

This optimization problem is non-convex, and there is no closed-form solution to it [5]. The brute force approach might be used to solve this problem. However, it is not practical to use the brute force approach once the number of parameters that need to be updated gets large.

It can be easily seen that the L0-norm minimization problem might not have a unique solution. However, other physical constraints can be applied to limit the solutions set, such as enforcing non-negative solutions only [6].

### 3.3 L1-norm minimization

L1-norm minimization has been used as a surrogate for L0-norm minimization [5]:

$$\min \|x\|_1 \ \ s.t \ Ax = b \qquad (4)$$

where the $\|x\|_1$ is the L1-norm of the vector $x$ :

$$\|x\|_1 = \sum_{i=1}^{n} |x_i| \tag{5}$$

The optimization problem in Eq. (4) is convex [5].

The noise and the error in the measurements can be taken as $e$ added to the right-side of Eq. (4), the alternative formulation shown in Eq. (6) can be used and modified with the free parameters [7].

$$\min \|x\|_1 \ s.t \ \|Ax - b\|_2 < \varepsilon \tag{6}$$

where $\varepsilon$ is a user-defined parameter to consider for the error in Eq. (6). $\varepsilon$ needs to be selected so that $\|e\|_2 < \varepsilon$. This means that the user needs to select a boundary on the error for the error for a stable recovery from the noisy estimated modal data.

L1-magic optimization toolbox [8] and CVX package [9] were used in this research for L1-norm minimization.

### 3.4 Lp-norm Minimization

The Lp-norm for a vector $x$ is defined as:

$$\|x\|_p = \left(\sum_{i=1}^{n} |x_i|^p\right)^{1/p} \tag{7}$$

where $0 < p < 1$. The minimization problem is:

$$\min \|x\|_p \ s.t \ Ax = b \tag{8}$$

Lp-norm minimization is a non-convex minimization problem that needs iterative algorithms, such as Iteratively Reweighted Least Squares (IRLS), to be solved [10]. It is claimed that the Lp-norm minimization provides a sparser solution and shows more robustness to noise compared to L1-norm minimization [10]. In our paper, Lp-norm minimization will be implemented in a damage characterization problem to compare its performance with L1-norm minimization algorithms, which have been used in this context before.

### 3.5 Effect of noise

In the previous research on using L1-norm minimization for damage characterization, the noise in the estimated modal features was very low [3]. The reason is that the sensitivity of modal features to damage is usually low. If the disturbance in the estimated modal features is bigger than the change in the modal features due to damage, it is not possible to achieve reasonable results.

### 3.6 Iterative methods

In order to find the solution to Eq. (1), linearization and iterative optimizations are required for a possible convergence [1]. Iterative methods have been used widely in the sensitivity method. However, they have not been investigated for undetermined systems in damage characterization. After linearizing Eq. (2), once $x$ (the change in the parameters) is found using L1-norm or Lp-norm minimization techniques, the parameters can be updated. The initial nominal model can be updated with the updated structural parameters and therefore a new set of modal features can be calculated with the updated model. However, there still will be a difference between the observed and the new features that were calculated from the updated model. The same procedure can be used to update the structural parameters and the model repeatedly until convergence.

An interesting question is whether the containing set of non-zeros (damage locations in our research) can change by using iterations. In all the investigated examples in this research, the set does not change. The iterations only adjust the non-zero entries. This can be helpful from the point of view of improving the solution to find the extent of the damage, but it also means that the iterations will not give a clue of other possible non-zero entries in the solution.

## 4 Numerical Illustrations

An example with the truss shown in Fig. 1 is considered. In this example, nominal modal is referred to model that needs to be updated and truth model is referred to the model that is used for simulation of the experiment. The truss is made of 20 bars and the bars are made of aluminum (The modulus of elasticity is 70 GPa). In the nominal model, the area of the cross-section for all bars is 0.01 $m^2$.

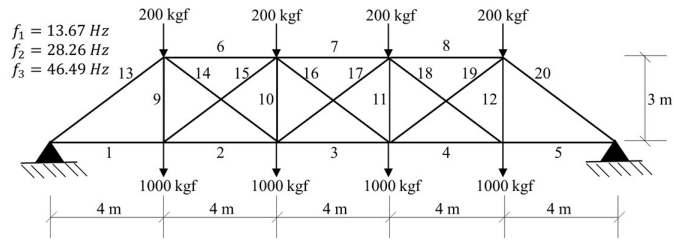

**Fig. 1.** Truss example used for damage characterization

In the truth model, the areas of bars 2 and 18 are reduced by 20% (which will reduce the bar stiffnesses by the same amount). The objective is to find the local damage in the bars by updating all the bar stiffnesses in the nominal model using only the frequencies. It is essential to mention that the difference in all the 16 frequencies between the truth model and nominal model is less than 3% (the average is 0.9%). This means the sensitivity of the frequencies with respect to the bar stiffnesses is very low. Therefore, if the estimation error in the frequencies is more than the change of the frequencies due to stiffness change, the damage cannot be identified regardless of the method that is used.

Fig. 2 and Fig. 3 show the results using L1-norm minimization with eight and nine estimated errorless frequencies. We also used iterations, as it was explained in section 3.6, to show its effect on the solution. The solution converged after four or five iterations. As it can be seen in these two figures, nine frequencies are needed to obtain a consistent solution with the actual physical damage. Using the eight frequencies or less will give a solution that is different than the physical solution. It is essential to mention that L1-norm minimization will find the solution with minimum L1-norm using any number of frequencies. However, this does not mean that the solution must be aligned with the physical damage. This means that the solution in Fig. 2 has a L1-norm that is smaller than the L1-norm of the physical solution.

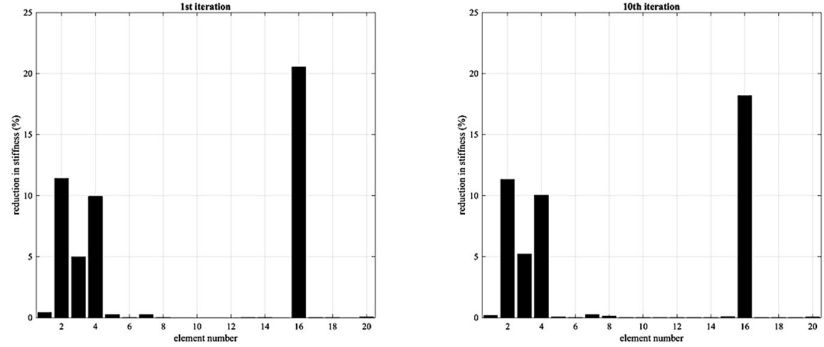

**Fig. 2.** Damage found in each element using eight frequencies

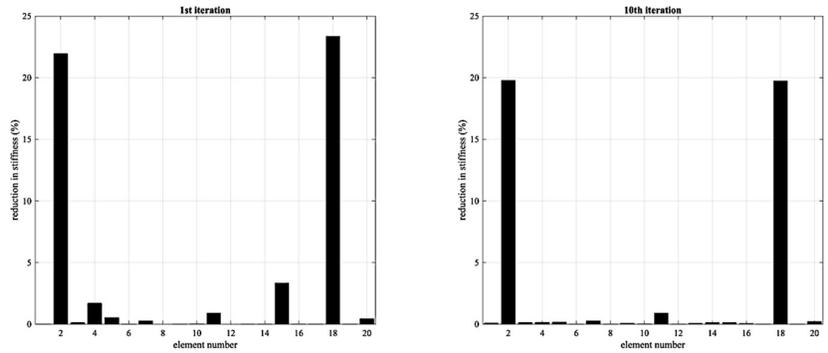

**Fig.3.** Damage found in each element using nine frequencies

As it can be seen in Fig. 2 and Fig. 3, the iterations did not change the location of non-zero entries in the solution, but they did improve the solution in terms of finding the extent of damage when the initial solution contained the true damaged bars.

In order to see the results for the case where the extent of damage is larger and therefore the nonlinear terms have more effects, 50% of damage is considered in bars 2 and 18. Fig. 4 shows the results using nine estimated frequencies.

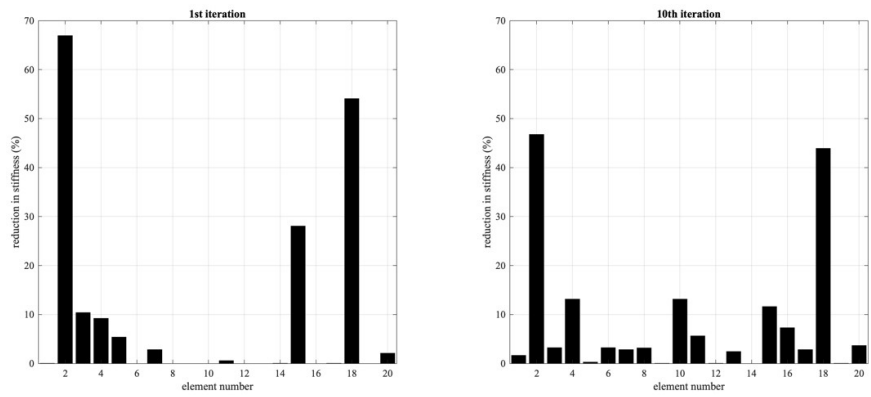

**Fig.4.** Damage found in each element using nine frequencies – Damage extent: 50%

As it can be seen, the results are affected by the nonlinear terms, and they are less impressive, but the iterations helped in terms of bringing the values of the damaged bars closer to the actual damage and reduce the extent of identified damage in bar 15 (which is zero). Based on these investigations, we can recommend using iterations in the sense that it will not hurt the results and can show some improvements in terms of finding the extent of the damage.

We also implemented Lp-norm (p=0.5) minimization using the IRLS method. In order to compare the performance of Lp-norm minimization with L1-norm minimization, we used 9 to 12 estimated frequencies under 5 different levels of noise to find the damage for the same scenario. The noise in the estimated frequencies is defined as:

$$no = \frac{f_m - f_e}{f_e} * 100 \qquad (5)$$

where $f_m$ is the noisy estimated frequency and $f_e$ is the exact errorless frequency of the truth model. Fig. 4 shows the success rate for 1000 realization for 1,2,3,4 and 5 percent uniform noise in the frequencies. In each realization, the value of noise for each frequency was drawn from a uniform distribution from 0 to $\pm no$.

In order to define "success" in identifying damage, a metric is needed. The primary objective in this example is the localization of the damage and then quantification. We call a solution, a successful solution, if non-zero elements in the solution contain the actual damage bars, and the number of non-zero elements in the solution is less than the number of estimated frequencies. If the solution has these two criteria in terms of localization of the damage, then it is possible to use iterations to enhance the quantification of the damage.

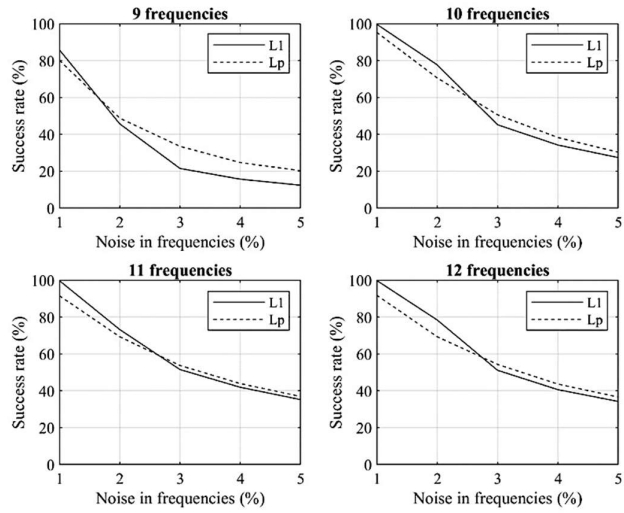

**Fig 5:** Damage characterization success rate in different levels of noise using different number of frequencies

As shown in Fig. 5, there is not a clear advantage in using Lp-norm minimization over L1-norm minimization. It might seem that Lp-norm minimization has a better performance in higher noise levels, but it should be mentioned that if the noise is more than 3 percent, both methods perform poorly.

## 5 Conclusions

In this study, we investigated the use of regularization techniques to exploit spareness in under-determined damage characterization problems. In previous research, the L1-norm minimization technique was used for regularization. However, the issue of noise and the need for using iterations for linearization were not highlighted in the previous research. Iterations have been applied and showed improvements in terms of finding the extent of the damage but they cannot usually help to localize the damage once the initial solution does not contain the damaged locations.

Lp-norm ($0<p<1$) minimization using IRLS was also implemented in this research which has not been used for damage characterization before. We compared the performance of Lp-norm minimization with L1-norm minimization to characterize damage in a truss example with the same level of noise in estimated frequencies. The results showed no clear advantage of using Lp-norm minimization over L1-norm minimization but it should be mentioned that using IRLS can be easier since there is no need estimate a boundary for the error.